\documentclass{amsart}
\usepackage{hyperref}
\usepackage{fullpage}
\usepackage{amsfonts}
\usepackage{amsthm}
\usepackage{amssymb}
\usepackage{amsmath}
\usepackage{mathtools}
\usepackage{tikz}
\usepackage{tikz-cd}
\usepackage{float}
\usepackage{mathrsfs} 
\usepackage{pdfpages}

\theoremstyle{definition}
\newtheorem{Th}{Theorem}[section]

\newtheorem{Prop}[Th]{Proposition}

\newtheorem{Rem}[Th]{Remark}
\newtheorem{?}[Th]{Problem}
\newtheorem{Lem}[Th]{Lemma}

\newtheorem*{PropP}{Proposition P}

\newcommand{\ovl}{\overline}
\newcommand{\wtd}{\widetilde}
\newcommand{\tb}{\textbf}
\newcommand{\ti}{\textit}
\newcommand{\Q}{\mathbb{Q}}
\newcommand{\Z}{\mathbb{Z}}
\newcommand{\C}{\mathbb{C}}
\newcommand{\R}{\mathbb{R}}

\newcommand{\A}{\mathbb{A}}

\newcommand{\ra}{\rightarrow}

\newcommand{\lan}{\langle}
\newcommand{\ran}{\rangle}
\newcommand{\mcl}{\mathcal}

\newcommand{\ts}{\otimes}
\newcommand{\be}{\begin{equation*}}
\newcommand{\ee}{\end{equation*}}
\newcommand{\mfk}{\mathfrak}
\newcommand{\bp}{\begin{proof}{\vspace{-\topsep}}}
\newcommand{\ep}{\end{proof}}
\newcommand{\e}{\epsilon}

\newcommand{\cond}{\text{cond}}
\newcommand{\udl}{\underline}
\def\len{8pt}

\DeclareMathOperator{\Gal}{Gal}
\DeclareMathOperator{\Tr}{Tr}
\DeclareMathOperator{\End}{End}
\DeclareMathOperator{\gal}{Gal}
\DeclareMathOperator{\sgn}{sgn}

\title{{Construction of Hecke Characters for Three-dimensional CM Abelian Varieties}}
\author{Zhengyuan Shang}
\date{\today}

\begin{document}

\maketitle
\begin{abstract}
It is well-known for an elliptic curve with complex multiplication that the existence of a $\Q$-rational model is equivalent to its field of moduli being equal to $\Q$, or its endomorphism ring being the ring of integers of 9 possible fields ($\ast$). Murabayashi and Umegaki proved analogous results for abelian surfaces. For three dimensional CM abelian varieties with rational fields of moduli, Chun narrowed down to a list of 37 possible CM fields. In this paper, we show that his list is exact. By constructing certain Hecke characters that satisfy a theorem of Shimura, we prove that precisely 28 isogeny classes of these abelian varieties have $\Q$-models. Therefore the complete analogy to $(\ast)$ fails here.
\end{abstract}

\section{Introduction}
For an elliptic curve $E$ defined over $\C$ with its endomorphism ring isomorphic to the ring of integers of an imaginary quadratic field $K$, it is a classical result that the followings are equivalent: \tb{i)} the $j$-invariant of $E$ is contained in $\Q$ (or equivalently, the field of moduli of $E$ is $\Q$), \tb{ii)} the class number $h_K$ of $K$ is one, \tb{iii)} $E$ has a model defined over $\Q$. The classification of class number one imaginary quadratic fields then enables us to determine all the (nine) $\Q$-rational points in the moduli space $\mcl{A}_1$ represented by CM elliptic curves. 

\par
Much effort has gone into generalizing this result to $n$-dimensional abelian varieties of CM-type for $n\geq 2$. In \cite{M1}, Murabayashi proved an analogue of \tb{i)} $\Leftrightarrow$ \tb{ii)} for polarized abelian surfaces: precise criterions were formulated for a surface to be simple and for its field of moduli to be $\Q$. Then in \cite{MU} and \cite{U}, all $\Q$-rational CM points in $\mcl{A}_2(d)$ for $d\geq 1$ were determined. Note that the direction $\tb{iii)}\Rightarrow \tb{i)}$ is clear by the definition of the field of moduli. By constructing certain Hecke characters over the quartic field, Murabayashi \cite{M2} showed that all such CM points have $\Q$-rational models, finishing the analogue to \tb{ii)} $\Rightarrow$ \tb{iii)}. 

Chun \cite{DC} gave necessary and sufficient conditions (Theorems \ref{thm: D1} \& \ref{thm: D2}) for three dimensional CM abelian varieties to have rational fields of moduli, and narrowed down to 37 potential fields  in Tables \ref{tab: T1} through \ref{tab: T4}. 

In this paper, we first show that these 37 fields all satisfy Chun's sufficient condition, so the abelian varieties corresponding to them admit polarizations such that their fields of moduli coincide with $\Q$. Then we prove that exactly 28 isogeny classes of them have $\Q$-models, indicating the failure of an analogue to \tb{i)} $\Leftrightarrow$ \tb{iii)}.

\par
\tb{Notation:} Unless otherwise specified, $K$ is a CM field of degree six over $\Q$ such that $G:=\gal(K/\Q)$ is cyclic with generator a $\sigma$, and $\rho=\sigma^3$ is complex conjugation. Let $F, F_0$ be the cubic and quadratic subfields of $K$. Let $U=\prod_u\mcl{O}_{K_u}^\times$ and $V=\prod_v\mcl{O}_{F_v}^\times$, where $u$ (resp. $v$) runs through all finite places of $K$ (resp. $F$). For a number field $k$, we use $\mcl{O}_k$, $E_k$, $W_k$, $C_k$, $I_k$, $H_k$, $P_k$, $h_k$, $f_k$, $d_k$, $k_{\A}^\times$, $k_\infty^\times$, $\mfk{f}$ to denote its ring of integers, group of units (in $\mcl{O}_k$), group of roots of unity, ideal class group, ideal group, Hilbert class field, defining polynomial (over $\Q$), class number, conductor, discriminant, idele group, Archimedian part (of $k_\A^\times$), and conductor over $F$. For a place $v$ of $k$, $k_v$ is the completion. For a fractional ideal $I\in I_k$, $cl(I)$, $N(I)$ are its class of in $C_k$ and absolute norm resp. For $x\in k_\A^\times$, $il(x)$ is the fractional ideal (in $k$) associated to $x$.
\be
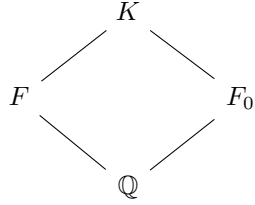
\begin{tikzcd}
& K & \\
F\ar[dash]{ru} \ar[dash]{rd} & & F_0 \ar[dash]{lu} \ar[dash]{ld}\\
& \Q &
\end{tikzcd}
\ee
Consider a structure $P=(A, C, \theta)$, where $A$ is a three dimensional abelian variety over $\C$ with a polarization $C$ and an injection $\theta: K\ra \End^0(A):=\End(A)\ts \Q$ such that $\theta^{-1}(\End(A))=\mcl{O}_K$. We always assume $\theta(K)$ is stable under the involution of $\End^0(A)$ determined by $C$. Let $\Phi=\{\sigma_1, \sigma_2, \sigma_3\}$ be a CM type of $K$ induced by the representation of K, through $\theta$, on the Lie algebra of $A$. Consider the isomorphism $\tilde{\Phi}: K\ts_\Q\R\ra\C^3$ such that $x\ts r \mapsto(rx^{\sigma_1}, rx^{\sigma_2}, rx^{\sigma_3})$. For $x\in K$, let $S_\Phi(x)$ be the diagonal matrix with $x^{\sigma_1}, x^{\sigma_2}, x^{\sigma_3}$ along the diagonal. By the theory of complex multiplication, there exists a fractional ideal $J$ of $K$ and an analytic isomorphism $\lambda: \C^3/\tilde{\Phi}(J)\ra A$ such that $\forall x\in \mcl{O}_K$, the following diagram commutes
\be
\begin{tikzcd}
\C^3/\tilde{\Phi}(J)\ar[r, "\lambda" ] \ar[d, swap, "S_{\Phi}(x)"]& A \ar[d, "\theta(x)"] \\
\C^3/\tilde{\Phi}(J)\ar[r, swap, "\lambda"] & A
\end{tikzcd}
\ee
Take a basic polar divisor in $C$ and consider its Riemann form $E(x, y)$ on $\C^3$ with respect to $\lambda$. Then there is an element $\eta$ of $K$ such that for $x, y\in K$, $E(\tilde{\Phi}(x), \tilde{\Phi}(y))=\Tr_{K/\Q}(\eta xy^\rho)$, $\eta^\rho=-\eta$, and \text{Im}$(\eta^{\sigma_j})>0$ for $1\leq j\leq 3$. We say $P$ is of type $(K, \Phi; \eta, J)$. Note that given $(K, \Phi; \eta, J)$ satisfying these conditions, we can also reconstruct the corresponding $(A, \mcl{C}, \theta)$.

\par
\tb{Acknowledgment:} I want to thank my mentor Prof. Matthias Flach for introducing me to this fascinating problem and offering numerous helpful suggestions. This project would not be possible without his guidance. I also wish to thank Caltech Class of 1936 for their generous support through the SURF program.

\section{Preliminaries}

In this section we state explicitly the results mentioned in the introduction, starting with Chun's criterions.

\Th{\cite[Theorem 1]{DC} \label{thm: D1} Let $P=(A, C, \theta)$ be of type $(K, \Phi; \eta, J)$. If $A$ is simple and the field of moduli of $(A, C)$ is $\Q$, then $K$ is one of the fields in Tables \ref{tab: T1} through \ref{tab: T4}, and $\Phi=\{1, \sigma, \sigma^2\}$ up to $G$-action.

\Rem{It is worth noting that fields} \label{thm: rm1} in Tables \ref{tab: T1} through \ref{tab: T4} are indeed cyclic over $\Q$. We denote by $A_i$ the field in Table \ref{tab: T1} with generator $\alpha_i$, and similarly the others by $B_i, \Gamma_i,$ and $\Omega_i$. By degree consideration, $|W_K|=2, 4, 6, 14,$ or 18. Moreover, $A_1=\Q(\zeta_{7})$ (resp. $A_2=\Q(\zeta_{9})$) is the only field with $|W_k|=14$ (resp. 18). Here $\zeta_l:=e^{\frac{2\pi i}{l}}$ for $l\in \Z^+$. By the uniqueness of the quadratic subfield  of $K$, for fields containing $i:=\sqrt{-1}$ (resp. $\omega:=e^{\frac{2}{3}\pi i}$), $|W_K|=4$ (resp. 6). For the rest, we then have $|W_K|=2$ and $W_K=\{1, -1\}$.

\Th{\cite[Theorem 2]{DC} \label{thm: D2} Let $P=(A, C, \theta)$ be of type $(K, \Phi; \eta, J)$}. If one of the followings is satisfied: \tb{i)} $K$ belongs to Table \ref{tab: T1}, \tb{ii)} $K$ belongs to Table \ref{tab: T2} and $H_{F}$ has four index 3 subfields, \tb{iii)} $K$ belongs to Table \ref{tab: T3} and $J$ is principal, \tb{iv)} $K=\Omega_1$, $cl(J)\in \ovl{C_{F}}$ (see the remark below), and $H_{F}$ has four index 3 subfields, and if $\Phi=\{1, \sigma, \sigma^2\}$ up to $G$-action, then $A$ is simple and the field of moduli of $(A, C)$ coincides with $\Q$. 

\Rem{By \cite[Lemma 2]{DC}, the canonical map $C_{F}\ra C_K$ is injective. Here $\ovl{C_{F}}$ is the image of $C_{F}$.

Moreover, we need the following theorem of Shimura, which is also crucial in Murabayashi's proof in \cite{M2}. Let $K\subset \C$ be a degree $2n$ CM field with $F$ its maximal real subfield. Let $\Phi$ be a CM type of $K$ and $(K', \Phi')$  be the reflex of $(K, \Phi)$. For a number field $k$ containing $K'$, consider $f: k^\times \ra K^\times$ with $f(x)=\prod_{\tau\in \Phi'}N_{k/K'}(x)^\tau$, which can be extended to a continuous homomorphism $k_\A^\times\ra K_\A^\times$ (also denoted $f$). Suppose $K, K'$ have subfields $D, D'$ respectively such that: \tb{i)} $D\subset F$, $D'\subset F'$, \tb{ii)} $K, K'$ are cyclic over $D, D'$ resp., \tb{iii)} $K=DS$, where $S/\Q$ is generated $\Tr \Phi'(a)=\sum_{\tau\in \Phi'}a^\tau$ for all $a\in K'$, and \tb{iv)} there is an isomorphism $\Gal(K'/D')\ra \Gal(K/D)$ where $\sigma\mapsto[\sigma]$ such that $(\Tr\Phi(a))^\sigma =\Tr\Phi(a^{[\sigma]})$ for all $a\in K$. 

\Th{\cite[Theorem 5]{GS}} \label{thm: S1} Let $A$ be an $n$-dimensional abelian variety over $\C$ with a polorization $\mcl{C}$ and $\theta: \mcl{O}_K\xhookrightarrow{} \End(A)$ such that the type of $(A, \theta)$ is $(K, \Phi)$.  Let $k_0$ be an algebraic number field containing both $D'$ and the field of moduli of $(A, \mcl{C}, \theta|_D)$. Let $k=k_0K'$ and $h=k_0F'$. Suppose that $K'$ and $k_0$ are linearly disjoint over $D'$ and let $\chi: h_\A^\times \ra\{\pm 1\}$ be the quadratic character associated to the extension $k/h$. Suppose there exists a Hecke character  $\psi: k_\A^\times \ra \C^\times$ satisfying
\begin{itemize}

\item[]
 \tb{a)} $x\in k_\infty^\times\Rightarrow \psi(x)=(f(x)_{\infty 1})^{-1}$, with $\infty 1$  the archimedian place of $K$ corresponding to $K\xhookrightarrow{} \C$
 \item[]
  \tb{b)} $y\in h_\A^\times\Rightarrow \psi(y)=\chi(y)|y|_{h_\A}^{-1}$, with $|\cdot|_{h_\A}$ the volume of ideles
  \item[]
  \tb{c)} $x\in k_\A^\times, x_\infty=1\Rightarrow \psi(x)\in K^\times, \psi(x)\psi(x)^\rho=|x|_{k_\A}^{-1}, (\psi(x))=il(f(x))$
  \item[]
   \tb{d)} $x\in k_\A^\times, x_\infty=1, \sigma\in \Gal(k/k_0)\Rightarrow \psi(x)^{[\sigma]}=\psi(x^\sigma)$, with $\Gal(k/k_0)$ and $\Gal(K'/D')$ identified.

\end{itemize} 
 Then there exists a structure $(A', \mcl{C}', \theta')$ rational over $k$ and isomorphic to $(A, \mcl{C}, \theta)$ such that corresponding Hecke character is $\psi$, and $(A', \mcl{C}', \theta'|_D)$ is rational over $k_0$.
 
\Rem{\label{thm: GS2} Let $\psi^\ast$ be the Hecke ideal character associated with $\psi$ and $\mfk{c}$ be the conductor of $\psi$. Then by \cite[Remark 2]{GS}, the conditions \tb{a)}, \tb{b}, \tb{c}, \tb{d)} are equivalent to 
\begin{itemize}
\item[] 
\tb{a$'$)} $\psi^\ast(s\mcl{O}_K)=\det \Phi'(N_{k/K'}(s))$, if $s\in k^\times$ and $s\equiv 1$ (mod $\mfk{c}$).

\item[]
\tb{b$'$)} $\psi^\ast(\mfk{t}\mcl{O}_K)=\left( \frac{k/h}{\mfk{t}}\right)N(\mfk{t})$, for every ideal $\mfk{t}$ of $h$ prime to $\mfk{c}$.

\item[]
\tb{c$'$)} For every ideal $I$ of $k$ prime to $\mfk{c}$,  $\psi^\ast(I)\in K^\times$, $\psi^\ast(I)\psi^\ast(I)^\rho=N(I)$, and $\psi^\ast(I)\mcl{O}_K=N_{k/K'}(I)^{\Phi'}$.

\item[]
\tb{d$'$)} For every ideal $I$ of $k$ prime to $\mfk{c}$ and $\sigma\in \Gal(k/k_0)$, $\psi^\ast(I)^{[\sigma]}=\psi^\ast(I^\sigma)$.
\end{itemize}
respectively. Here the equivalence \tb{c)}$\Leftrightarrow\tb{c$'$)}$ requires that the lattice $J$ in $K$ determined by $A$ is a fractional ideal, which holds in our case as $\End(A)= \theta(\mcl{O}_K)$.

Finally, we have the following converse to Theorem \ref{thm: S1} by Yoshida.  

\Th{\label{thm: S2} \cite[Theorem 2, Corollary]{HY}} Fix a structure $(A, \mcl{C}, \theta)$ as above. If in the isogeny class of $(A, \theta)$, there exists a structure $(A', \theta')$ such that $(A', \theta'|_{D})$ is rational over $k_0$ and $\End(A')=\theta'(\mcl{O}_K)$, then there exists a Hecke character $\psi$ of $k_\A^\times$ that satisfies \tb{a)} through \tb{d)}.

\section{Main Results}

\Th{Let $(A, \theta)$ be of type $(K, \Phi)$. Then $A$ is simple and the field of moduli of $(A, \mcl{C})$ is $\Q$ for some polarization $\mcl{C}$ if and only if $K$ is one of the fields in Tables \ref{tab: T1} through \ref{tab: T4}, and $\Phi=\{1, \sigma, \sigma^2\}$ up to $G$-action.}

\bp
The $\Rightarrow$ direction is a restatement of Theorem \ref{thm: D1}. For $F$ in Tables \ref{tab: T2} and \ref{tab: T4}, we have $h_{F}=3$, so $H_{F}$ is Galois over $\Q$ of degree 9. In particular, $\gal(H_F/\Q)$ is either $\Z/9\Z$ or $\Z/3\Z\times \Z/3\Z$. Since $H_F$ is abelian, it coincides with the genus field of $F$, which is the compositum of cubic fields by \cite[Theorem]{ISH}. Hence $\gal(H_F/\Q)\cong \Z/3\Z\times \Z/3\Z$ and the converse follows from Theorem \ref{thm: D2}.
\ep

\vspace{-\len}
Now let $K$ be one of the fields in Tables \ref{tab: T1} through \ref{tab: T4} and $\Phi=\{1, \sigma, \sigma^2\}$. Choose $\theta$ such that $(A, \theta)$ is of type $(K, \Phi)$. Then there exists a polarization $\mcl{C}$ such that the field of moduli of $(A, \mcl{C}, \theta)$ is $\Q$. Since $(K, \Phi)$ is simple, $(K', \Phi')=(K, \{1, \sigma^4, \sigma^5\})$. To apply Theorem \ref{thm: S1}, take $D=D'=\Q$. Note that $\Q\subset F = F'$ and $K=K'$ is cyclic over $\Q$. Meanwhile, the field over $\Q$ generated by $\sum_{\tau\in \Phi'}a^\tau$ for all $a\in K'$ is just $K$, and $K=K\Q$. For the identity isomorphism $\sigma\ra [\sigma]$ from $\gal(K'/D')\ra \gal(K/D)$, we have $\Tr(\Phi(a))^\sigma=\Tr \Phi(a^{[\sigma]})$, $\forall a\in K$. Put $k_0=\Q$, which contains both $D'$ and the field of moduli of $(A, \mcl{C}, \theta_D)$. Then $k=K$ and $h=F$. Note that $K'$ and $k_0$ are linearly disjoint over $D'$. Furthermore, $f: K^\times_\A\ra K^\times_\A$ is defined by $x\mapsto xx^{\sigma^4}x^{\sigma^5}$.

\PropP{There exists a Hecke character $\psi$ of $K_\A^\times\ra \C^\times $  satisfying the conditions \tb{a)}, \tb{b)}, \tb{c)}, \tb{d)} if and only if $|W_K|\neq 6$ (or equivalently, $K\notin\{A_4, A_7, A_{15}, A_{17}, B_1, B_2, B_5, \Gamma_5, \Gamma_6\})$. 

\bp
This follows from Propositions \ref{thm: cyc}, \ref{thm: p1}, \ref{thm: p2}, \ref{thm: p3}, \ref{thm: p4}, and \ref{thm: no}.
\ep

\Rem{In our case, \tb{a)}, \tb{b)}, \tb{c)}, \tb{d)} translate to}

\begin{itemize}

\item[]
 \tb{a$''$)} $x\in K_\infty^\times\Rightarrow \psi(x)=(xx^{\sigma^4}x^{\sigma^5})_{\infty 1}^{-1}$
 \item[]
  \tb{b$''$)} $y\in F_\A^\times\Rightarrow \psi(y)=\chi(y)|y|_{F_\A}^{-1}$
  \item[]
  \tb{c$''$)} $x\in K_\A^\times, x_\infty=1\Rightarrow \psi(x)\in K^\times, \psi(x)\psi(x)^\rho=|x|_{K_\A}^{-1}, (\psi(x))=il(xx^{\sigma^4}x^{\sigma^5} )$
  \item[]
   \tb{d$''$)} $x\in K_\A^\times, x_\infty=1, \sigma\in \Gal(K/\Q)\Rightarrow \psi(x)^{\sigma}=\psi(x^\sigma)$

\end{itemize}

\Th{Let $(A, \theta)$ be of type $(K, \Phi)$. If $K$ is one of the fields in Tables \ref{tab: T1} through \ref{tab: T4}, and $\Phi=\{1, \sigma, \sigma^2\}$ up to $G$-action, then $(A, \theta)$ has a $\Q$-model if and only if $|W_K|\neq 6$.}
\bp
This is immediate from Proposition \tb{P}, Theorem \ref{thm: S1}, and Theorem \ref{thm: S2}.
\ep

\Rem{If} $|W_K|\neq 6$, then the $\Q$-model of $(A, \theta)$ has conductors as listed in Tables \ref{tab: T1} through \ref{tab: T4} by Proposition \ref{thm: cond}. Otherwise, $(A, \theta)$ has a $F$-model by Proposition \ref{thm: last}.}

\section{Hecke Characters}

We prove Proposition \tb{P} in this section, starting with the two cyclotomic fields $K$ with $|W_K|=14$ or 18.

\Prop{\label{thm: cyc} Proposition \tb{P} holds for $K=\Q(\zeta_7), \Q(\zeta_9)$ (i.e. $A_1, A_2$).}

\bp

The case $K=A_1=\Q(\zeta_{7})$ follows immediately from Shimura's construction in \cite[Theorem 7]{GS}. 
If $K=A_2=\Q(\zeta_9)$, we adopt a similar approach. Consider the prime ideal $\mfk{l}=(1-\zeta_{9})\mcl{O}_K$ above $(3)$. Then $3\mcl{O}_K=\mfk{l}^6$. By the theory of complex multiplication (cf. \cite[Lemma 5]{DC}), for an ideal $I$ in $K$, $I^{\Phi'}=II^{\sigma^4}I^{\sigma^5}=r\mcl{O}_K$ for some $r\in K$ such that $rr^\rho=N(I)$. If $I$ is prime to $3$, then $N(I)\equiv 1$ (mod 3), so $r\equiv \pm 1$ (mod $\mfk{l}\mcl{O}_{K_{\mfk{l}}}$). Note that $1, \zeta_9, \zeta_9^2, ..., \zeta_9^8\in 1+\mfk{l}\mcl{O}_{K_{\mfk{l}}}$ project to distinct classes in $(1+\mfk{l}\mcl{O}_{K_{\mfk{l}}})/(1+\mfk{l}^4\mcl{O}_{K_{\mfk{l}}})$. Moreover, $1, \zeta_9, \zeta_9^2$ map to distinct classes in $(1+\mfk{l}\mcl{O}_{K_{\mfk{l}}})/(1+\mfk{l}^2\mcl{O}_{K_{\mfk{l}}})$ and $1, \zeta_9^3, \zeta_9^6$ to distinct classes in $(1+\mfk{l}^3\mcl{O}_{K_{\mfk{l}}})/(1+\mfk{l}^4\mcl{O}_{K_{\mfk{l}}})$.
If $r\equiv 1$ (mod $\mfk{l}^2\mcl{O}_{K_{\mfk{l}}}$), suppose $r=1+t(1-\zeta_9)^2$ for some $t\in \mcl{O}_{K_\mfk{l}}$. Then since $t\equiv t^\rho$ (mod $\mfk{l}\mcl{O}_{K_{\mfk{l}}}$), we have  $rr^\rho\equiv (1+t(1-\zeta_9)^2)(1+t^\rho(1-\zeta_9^8)^2)\equiv
 (1+t(1-\zeta_9)^2)(1+t(1-\zeta_9^2)^2)\equiv 1+t(2-\zeta_9-\zeta_9^2)\equiv 1-t(\zeta_9-1)^2\equiv 1$ (mod $\mfk{l}^3\mcl{O}_{K_{\mfk{l}}}$). Thus $t\equiv 0$ (mod $\mfk{l}\mcl{O}_{K_{\mfk{l}}}$) and $r\equiv 1$ (mod $\mfk{l}^3\mcl{O}_{K_{\mfk{l}}})$. Consequently, there is a unique sign of $\pm$ and an integer $0\leq m\leq 8$ such that $\pm \zeta_9^mr\equiv 1$ (mod $\mfk{l}^4\mcl{O}_{K_{\mfk{l}}}$). We then define $\psi^\ast(I)=\pm \zeta_9^mr$. Since $rr^\rho=N(I)$ and $\{\pm \zeta_9^m\}$ are all the roots of unity, $\psi^\ast$ is well-defined.

\par
Note that $\psi^\ast$ is a homomorphism from the group $I_{\mfk{l}^4}$ of ideals of $K$ prime to $\mfk{l}^4$ to $K^\ast$. For $(s)\in I_{\mfk{l}^4}$ with $s\equiv 1$ (mod $\mfk{l}^4$), since $\sigma$ fixes $\mfk{l}$, we have $ss^{\sigma^4}s^{\sigma^5}\equiv 1$ (mod $\mfk{l}^4$). Then $\psi^\ast((s))=ss^{\sigma^4}s^{\sigma^5}$, so $\psi^\ast$ is a Hecke character. By construction, it satisfies \tb{a$'$)}, \tb{c$'$)}, and \tb{d$'$)} in Remark \ref{thm: GS2}. For \tb{b$'$)}, it suffices to show that 
\be
\psi^\ast(\mfk{p}\mcl{O}_K)=\left( \frac{K/F}{\mfk{p}}\right) N(\mfk{p})
\tag{$*$}
\ee
for every prime ideal $\mfk{p}$ of $F$ prime to $3$. Such $\mfk{p}$ is unramified in $K$. Let $\rho_0$ be the generator of $\Gal(K/F)$. If $\mfk{p}\mcl{O}_K=\mfk{q}\mfk{q}^{\rho_0}$ for some prime ideal $\mfk{q}$ of $K$, then $\psi^\ast(\mfk{p}\mcl{O}_K)=\psi^\ast(\mfk{q})\psi^\ast(\mfk{q})^{\rho_0}=N(\mfk{q})=N(\mfk{p})$, so $(\ast)$ holds. Otherwise, $\mfk{p}$ remains prime in $K$. Then $(\frac{K/F}{\mfk{p}})=\rho_0$ and $N(\mfk{p})\equiv -1$ (mod 3). 
Meanwhile, $(\mfk{p}\mcl{O}_K)^{\Phi'}=N(\mfk{p})\mcl{O}_K$, so $\psi^\ast(\mfk{p}\mcl{O}_K)=-N(\mfk{p})$ and $(\ast)$ also holds. This completes the proof.
\ep

\vspace{-\len}
For the remaining fields from Tables 1 through 4, $|W_K|\in \{2, 4, 6\}$. We first show the existence of satisfactory Hecke characters $\psi$ for the cases when $|W_K|=2$ or $4$, where the following proposition is an essential step. 

\Prop{\label{thm: equiv} Suppose $|W_K|=2$ or 4. If there exists a character $\wtd{\chi}_{\mfk{f}}: \prod_{u|\mfk{f}}\mcl{O}_{K_u}^\times \ra W_K$ such that 
\begin{align*}
&\wtd{\chi_\mfk{f}}|_{\prod_{v|\mfk{f}}\mcl{O}_{F_v}^\times}=\chi|_{\prod_{v|\mfk{f}}\mcl{O}_{F_v}^\times}\\
&\wtd{\chi_\mfk{f}}(x^\sigma)= \wtd{\chi_\mfk{f}}(x)^\sigma\\
&\wtd{\chi_\mfk{f}}(x)=x, \forall x\in W_K
\end{align*}
where $u$ and $v$ are finite places of $K$ and $F$. Then there exists a continuous homomorphism $\psi_0: K^\times U K^\times_\infty\ra \C^\times$ satisfying
\tb{a$''$)}, \tb{b$''$)}, \tb{c$''$)}, \tb{d$''$)} with $K^\times_\A$, $F^\times_\A$ replaced by $K^\times UK^\times_\infty$, $F^\times VF^\times_\infty$ respectively.

\bp
Since $F$ is the unique cubic subfield of $K$, $\mfk{f}^\sigma=\mfk{f}$, so the Galois equivariance condition on $\wtd{\chi_\mfk{f}}$ is well-defined. We adopt the construction of Yoshida in \cite[Section 4]{HY}. For $x\in F_\A^\times$ with $x_v\in \mcl{O}_{F_v}^\times$ for each finite place $v $ of $F$ not dividing $\mfk{f}$, we have 
\be
\chi(x)=\prod_{v|\mfk{f}}\chi_v(x_v)\sgn(x_{\infty 1})x\sgn(x_{\infty 2})\sgn(x_{\infty 3})
\ee
where $\sgn: \R^\times \ra \{\pm 1\}$ is the sign map. Now consider $\psi_0: K^\times UK_\infty^\times\ra \C^\times$ defined by 
\be
\psi_0(xyz)=\wtd{\chi_\mfk{f}}(y_\mfk{f})(f(z))_{\infty 1}^{-1}
\ee
where $x\in K^\times, y\in U, z\in K^\times_\infty$, and $y_\mfk{f}$ is the $\prod_{u|\mfk{f}}\mcl{O}_{K_u}^\times$-component of $y$. Note that $K^\times \cap UK^\times_\infty=E_K=W_KE_F$ by \cite[Lemma 2]{DC}. For $x\in E_F$, $x'\in W_K$, we have 
\be
\wtd{\chi_\mfk{f}}(xx')(f(xx'))^{-1}_{\infty 1}= x'(f(x'))_{\infty 1}^{-1}\prod_{v|\mfk{f}}\chi_v(x)N_{F/\Q}(x)=\chi(x)=1
\ee
Hence $\psi_0$ is well-defined. Note that \tb{a$''$)} is satisfied. For $x\in F^\times, y\in V, z=(z_1, z_2, z_3)\in F^\times_\infty$, we have 
\be
\psi_0(xyz)=\wtd{\chi_\mfk{f}}(y_\mfk{f})(f(z))_{\infty 1}^{-1}=\prod_{v|\mfk{f}}\chi_v(y_v)(z_1z_2z_3)^{-1}
\ee
Meanwhile, 
\be
\chi(xyz)|xyz|_{F_\A}^{-1}=\prod_{v|\mfk{f}}\chi_v(y_v)\sgn(z_1)\sgn(z_2)\sgn(z_3)N(il(y))|z_1z_2z_3|^{-1}=\prod_{v|\mfk{f}}\chi_v(y_v)(z_1z_2z_3)^{-1}
\ee
so \tb{b$''$)} holds. For $x\in K^\times, y\in U, z=(z_1, z_2, z_3)\in K^\times_\infty$, if $(xyz)_\infty=1$, then $xz_1=x^\sigma z_2=x^{\sigma^2}z_3=1$, so 
\be
\psi_0(xyz)=\wtd{\chi_\mfk{f}}(y_\mfk{f})(f(z))_{\infty 1}^{-1}=\wtd{\chi_\mfk{f}}(y_\mfk{f})(z_1z_2^{\sigma^4}z_3^{\sigma^2})^{-1}=\wtd{\chi_\mfk{f}}(y_\mfk{f})xx^{\sigma^4}x^{\sigma^5}
\ee
Thus $\psi_0(xyz)\psi_0(xyz)^\rho=N_{K/\Q}(x)=N(il(xyz))=|xyz|^{-1}_{K_\A}$. Moreover, we have $(\psi_0(xyz))=(xx^{\sigma^4}x^{\sigma^5})=il(f(xyz))$, so \tb{c}$''$) holds. Note that 
\begin{align*}
\psi_0((xyz)^\sigma)&=\wtd{\chi_\mfk{f}}((y^\sigma)_\mfk{f})(f(z^\sigma))_{\infty 1}^{-1}=\wtd{\chi_\mfk{f}}((y_\mfk{f})^\sigma)(z_1z_2^{\sigma^4}z_3^{\sigma^5})^{-1}\\
&=\wtd{\chi_\mfk{f}}(y_\mfk{f})^\sigma xx^\sigma x^{\sigma^5}=(\wtd{\chi_\mfk{f}}(y_\mfk{f})xx^{\sigma^4}x^{\sigma^5})^\sigma=(\wtd{\chi_\mfk{f}}(y_\mfk{f})(z_1z_2^{\sigma^4}z_3^{\sigma^2})^{-1})^\sigma=\psi_0(xyz)^\sigma
\end{align*}
Hence \tb{d$''$)} also holds, which completes the proof.
\ep

\Lem{\label{thm: quad} For fields $K$ with $|W_K|=2$ and $F_0\neq \Q(\sqrt{-2})$, there exists $\wtd{\chi_\mfk{f}}$ satisfying the conditions above.}

\bp
Since the rational prime (2) is unramified in both $F$ and $F_0$, it is unramified in $K$ and $(2, \mfk{f})=1$. For a prime $\mfk{p}\subset F$ dividing $\mfk{f}$, let $\mfk{P}\subset K$ be the prime above. Note that $\chi_\mfk{p}|_{\mcl{O}_{F_\mfk{p}}^\times}$ is the quadratic symbol mod $\mfk{p}$ in $\mcl{O}_{F_\mfk{p}}^\times$. Let $\wtd{\chi_\mfk{P}}$ be the quadratic symbol mod $\mfk{P}$ in $\mcl{O}_{K_\mfk{P}}^\times$. Then $\wtd{\chi_\mfk{f}}=\prod_{\mfk{P}|\mfk{f}}\wtd{\chi_\mfk{P}}$ is the desired character.
\ep

\vspace{-\len}
To construct $\wtd{\chi_\mfk{f}}$ for fields with $F_0=\Q(\sqrt{-2})$ or $|W_K|=4$ (i.e. $F_0=\Q(i)$), we make the following observation.

\Lem{\label{thm: L} For $K=A_5, A_6, A_9, A_{12}, A_{16}$, the prime ideal $(2)$ is inert in $F$, while for $\Gamma_1, \Gamma_4, \Gamma_8, \Omega_1$, it splits completely. Let $\mfk{p}\subset F$ be a prime above $(2)$. Then $\mfk{p}\mcl{O}_K=\mfk{P}^2$ for some prime  $\mfk{P}\subset K$, and $\mfk{f}$ is a power of $\mfk{p}$ or $\mfk{p}_1\mfk{p}_2\mfk{p}_3$. Meanwhile, the quadratic character $\chi_\mfk{p}: F_\mfk{p}^\times \ra \{\pm 1\}$ induced by $K_\mfk{P}/F_{\mfk{p}}$ is trivial on $1+\mfk{p}^r$, where $r=2$ for $A_5, A_6, A_{12}, \Gamma_1, \Gamma_4, \Omega_1$ and $r=3$ for $A_9, A_{16}, \Gamma_8$ (corresponding to $F_0=\Q(i)$ or $\Q(\sqrt{-2})$).

\bp
The first two claims follow from \cite[Theorem 1]{PE}. Since $F$ and $F_0$ are cubic and quadratic respectively, $d_F=f_F^2$ and $d_{F_0}=f_{F_0}$ are coprime. As $K/F$ is quadratic, $\mfk{f}=d_{K/F}$, where $d_{K/F}$ is the relative discriminant. Therefore $d_F^2N_{F/\Q}(\mfk{f})=d_K=d_F^2d_{F_0}^3$ and $N_{F/\Q}(\mfk{f})=d_{F_0}^3$, which proves the last two claims.
\ep

\vspace{-\len}
In Lemmas \ref{thm: L1} through \ref{thm: L4}, we consider $\mfk{p}, \mfk{P}, \chi_\mfk{p}$ as defined above.

\Lem{\label{thm: L1} [$F_0=\Q(\sqrt{-2})$ and (2) splits in $F$] For $K=\Gamma_8$, there exists a homomorphism $\wtd{\chi_\mfk{P}}: \mcl{O}_{K_\mfk{P}}^\times\ra \{\pm 1\}$ such that}
$\wtd{\chi_\mfk{P}}|_{\mcl{O}_{F_\mfk{p}}^\times} = \chi_\mfk{p}|_{\mcl{O}_{F_\mfk{p}}^\times}$ and 
$\wtd{\chi_\mfk{P}}(x^\rho) = \wtd{\chi_\mfk{P}}(x)^\rho, \forall x\in \mcl{O}_{K_\mfk{P}^\times}$.

\bp
Note that $\chi_\mfk{p}$ (when restricted to ${\mcl{O}_{F_\mfk{p}}^\times}$) factors through $\mcl{O}_{F_\mfk{p}}^\times/(1+\mfk{p}^3)\cong (\mcl{O}_{F_\mfk{p}}/\mfk{p}^3)^\times$, so it induces a homomorphism $\ovl{\chi_\mfk{p}}: (\mcl{O}_{F_\mfk{p}}/\mfk{p}^3)^\times \ra \{\pm 1\}$. Since $\mcl{O}_{F_\mfk{p}}\cap \mfk{P}^5=\mfk{p}^3$, there is a canonical inclusion $(\mcl{O}_{F_\mfk{p}} /\mfk{p}^3)^\times \ra (\mcl{O}_{K_\mfk{P}}/\mfk{P}^5)^\times$. We construct a homomorphism $\ovl{\chi_\mfk{P}}: (\mcl{O}_{K_\mfk{P}}/\mfk{P}^5)^\times \ra \{\pm1\}$ lifting $\ovl{\chi_\mfk{p}}$ that is also compatible with the action of $\rho$. Let $\eta=\sqrt{-2}$, then $\eta$ is a uniformizer of $\mcl{O}_{K_\mfk{P}}$ and $\eta^\rho=-\eta$. Since $(2)$ splits completely, we can take $F_\mfk{p}=\Q_2$ and $K_\mfk{P}=\Q_2(\sqrt{-2})$. Let $\xi_1=\ovl{1+\eta}, \xi_2=\ovl{1+\eta^2}, \xi_3=\ovl{1+\eta^4}\in (\mcl{O}_{K_\mfk{P}}/\mfk{P}^5)^\times$. Then 
\begin{align*}
(\mcl{O}_{K_\mfk{P}}/\mfk{P}^5)^\times & =\lan \xi_1 \ran \times \lan \xi_2 \ran \times \lan \xi_3\ran  \cong \Z/4\Z \times \Z/2\Z \times \Z/2\Z \\
(\mcl{O}_{F_\mfk{p}}/\mfk{p}^3)^\times &= \lan \xi_2 \ran \times \lan \xi_3 \ran \cong \Z/2\Z\times \Z/2\Z 
\end{align*}
Let $\ovl{\chi_\mfk{P}}(\xi_1)=1$, $\ovl{\chi_\mfk{P}}(\xi_j)=\ovl{\chi_\mfk{p}}(\xi_j)$ for $j= 2, 3$.
As $N_{K_\mfk{P}/F_\mfk{p}}(5+5\eta)\equiv (1+\eta^2)(1+\eta^4)$ (mod $\mfk{p}^3$), $\ovl{\chi_\mfk{p}}(\xi_2\xi_3)=1$. Since $\xi_1^\rho=\xi_1^3\xi_2\xi_3$,  $\ovl{\chi_\mfk{P}}(\xi^\rho_1)=\ovl{\chi_\mfk{P}}(\xi_1)^\rho$. Then the composition of $\ovl{\chi_\mfk{P}}$ with the projection satisfies the claim.
\ep

\Lem{\label{thm: L2} [$F_0=\Q(\sqrt{-2})$ and $(2)$ is inert in $F$] For $K=A_9, A_{16}$, there exists a homomorphism $\wtd{\chi_\mfk{P}}: \mcl{O}_{K_\mfk{P}}^\times\ra \{\pm 1\}$ such that}
$\wtd{\chi_\mfk{P}}|_{\mcl{O}_{F_\mfk{p}}^\times} = \chi_\mfk{p}|_{\mcl{O}_{F_\mfk{p}}^\times}$ and
$\wtd{\chi_\mfk{P}}(x^\sigma) = \wtd{\chi_\mfk{P}}(x)^\sigma, \forall x\in \mcl{O}_{K_\mfk{P}^\times}$.

\bp
Note that $F_\mfk{p}$ is the unramied cubic extension of $\Q_2$ and $K_\mfk{P}$ is the composition of $F_\mfk{p}$ and $\Q_2(\sqrt{-2})$. Thus we have the following commutative diagram, where the horizontal arrows are local Artin maps

\be
\begin{tikzcd}
F_\mfk{p}^\times \ar[r, "\chi_\mfk{p}"] \ar[d, swap, "N_{F_\mfk{p}/\Q_2}"]& \gal(K_\mfk{P}/F_\mfk{p}) \ar[d] \\
\Q_2^\times \ar[r] \ar[rd]& \gal(K_\mfk{P}/\Q_2) \ar[d] \\
 & \gal(\Q_2(\sqrt{-2})/\Q_2)
\end{tikzcd}
\ee
Let $\wtd{\chi_\mfk{P}}$ be the composition of $N_{K_\mfk{P}/\Q_2(\sqrt{-2})}$ with the homomorphism $\mcl{O}_{\Q_2(\sqrt{-2})}^\times \ra \{\pm 1\}$ defined in Lemma \ref{thm: L1}. Then $\wtd{\chi_\mfk{P}}$ satisfies our requirement.
\ep

\Lem{\label{thm: L3} [$F_0=\Q(i)$ and $(2)$ splits in $F$] For $K= \Gamma_1, \Gamma_4, \Omega_1$, there exists a homomorphism $\wtd{\chi_\mfk{P}}: \mcl{O}_{K_\mfk{P}}^\times\ra \{\pm 1, \pm i\}$ such that 
$\wtd{\chi_\mfk{P}}|_{\mcl{O}_{F_\mfk{p}}^\times}=\chi_\mfk{p}|_{\mcl{O}_{F_\mfk{p}}^\times}$, $\wtd{\chi_\mfk{P}}(i)=i$, and $\wtd{\chi_\mfk{P}}(x^\rho)=\wtd{\chi_\mfk{P}}(x)^\rho, \forall x\in \mcl{O}_{K_\mfk{P}^\times}$.

\bp 
By the reasoning in Lemma \ref{thm: L1}, it suffices to contruct a suitable lift $\ovl{\chi_\mfk{P}}$ of $\ovl{\chi_\mfk{p}}$ for $K=\Gamma_1$. We can take $F_\mfk{p}, K_\mfk{P}$ as $\Q_2$ and $\Q_2(i)$ resp. Let $\eta=1-i$. Consider $\xi_1=\ovl{1+\eta}, \xi_2=\ovl{1+\eta^2}\in (\mcl{O}_{K_\mfk{P}}/\mfk{P}^3)^\times$. Then
\begin{align*}
(\mcl{O}_{K_\mfk{P}}/\mfk{P}^3)^\times & = \lan \xi_1\ran \times \lan \xi_2\ran \cong \Z/2\Z \times \Z/2\Z \\
(\mcl{O}_{F_\mfk{p}}/\mfk{p}^2)^\times & = \lan \xi_2 \ran \cong \Z/2\Z
\end{align*}
Let $\ovl{\chi_\mfk{P}}(\xi_1)=i$ and $\ovl{\chi_\mfk{P}}(\xi_2)=\ovl{\chi_\mfk{p}}(\xi_2)=-1$, so $\ovl{\chi_\mfk{P}}$ lifts $\ovl{\chi_\mfk{p}}$. Since $\xi_1^\rho=\xi_1\xi_2$, $\ovl{\chi_\mfk{P}}(\xi_1^\rho)=\ovl{\chi_\mfk{P}}(\xi_1)^\rho$. Moreover, $\wtd{\chi_\mfk{P}}(i)=\ovl{\chi_\mfk{P}}(\xi_1)=i$. This completes the proof.
\ep

\Lem{\label{thm: L4} [$F_0=\Q(i)$ and $(2)$ is inert in $F$] For $K= A_5, A_6, A_{12}$, there exists a homomorphism $\wtd{\chi_\mfk{P}}: \mcl{O}_{K_\mfk{P}}^\times\ra \{\pm 1, \pm i\}$ such that 
$\wtd{\chi_\mfk{P}}|_{\mcl{O}_{F_\mfk{p}}^\times}=\chi_\mfk{p}|_{\mcl{O}_{F_\mfk{p}}^\times}$, $\wtd{\chi_\mfk{P}}(i)=i$, and $\wtd{\chi_\mfk{P}}(x^\sigma)=\wtd{\chi_\mfk{P}}(x)^\sigma, \forall x\in \mcl{O}_{K_\mfk{P}^\times}$. 

\bp
Let $\wtd{\chi_\mfk{P}}$ be the composition of $N_{K_\mfk{P}/\Q_2(i)}\circ \rho$ with the homomorphism $\mcl{O}_{\Q_2(i)}^\times \ra \{\pm 1, \pm i\}$ above. 
\ep

\Rem{\label{thm: sys} When $K=\Gamma_1, \Gamma_4, \Omega_1$, suppose $2=\mfk{p}_1\mfk{p}_2\mfk{p}_3=\mfk{P}_1^2\mfk{P}_2^2\mfk{P}_3^2$ for prime ideals $\mfk{p}_j\subset F$ and $\mfk{P}_j\subset K$ such that $\sigma\mfk{P}_1=\mfk{P}_2$ and $\sigma\mfk{P}_2=\mfk{P}_3$.} Then $\sigma$ induces isomorphisms $K_{\mfk{P}_1}\ra K_{\mfk{P}_2}$ and $K_{\mfk{P}_2}\ra K_{\mfk{P}_3}$. Let $\varphi: K_{\mfk{P}_3}\ra \Q_2(i)$ be an isomorphism such that $\varphi(i)=i$. 
\be
\begin{tikzcd}
 & \Q_2(i) & \\
 K_{\mfk{P}_1} \ar[r, swap, "\sigma"] \ar[ur] & K_ {\mfk{P}_2} \ar[r, swap, "\sigma"] \ar[u]& K_{\mfk{P}_3} \ar[ul, swap, "\varphi"]
\end{tikzcd}
\ee
Let $\udl{\chi}$ be the character of $\mcl{O}_{\Q_2(i)}^\times$ defined in Lemma \ref{thm: L3}. Consider $\wtd{\chi_\mfk{f}}: \prod_{j=1}^3\mcl{O}_{K_{\mfk{P}_j}}^\times \ra \{\pm 1, \pm i\}$ with
\be
\wtd{\chi_\mfk{f}}(y_1, y_2, y_3) = \udl{\chi}(\varphi(y_1^{\sigma^5}))\udl{\chi}(\varphi(y_2^{\sigma^1}))\udl{\chi}(\varphi(y_3^\rho))
\ee
Then $\wtd{\chi_\mfk{f}}$ restricts to $\prod_{j=1}^3\chi_{\mfk{p}_j}$ on $\prod_{j=1}^3\mcl{O}^\times_{F_{\mfk{p}_j}}$, $\wtd{\chi_\mfk{f}}(i, i, i)=\udl{\chi}(-i)^3=i$ and $\wtd{\chi_\mfk{f}}((y_1, y_2, y_3)^\sigma)=\wtd{\chi_\mfk{f}}(y_3^\sigma, y_1^\sigma, y_2^\sigma)=\wtd{\chi_\mfk{f}}(y_1, y_2, y_3)^\sigma$. Note that this is simply the composition $\udl{\chi}\circ N_{K/\Q(i)}\circ \rho$ on $K\ts \Q_2$. Similarly, when $K=\Gamma_8$, we can construct a $\wtd{\chi_\mfk{f}}$ from the character  in Lemma \ref{thm: L1} that satisfies the conditions of Proposition \ref{thm: equiv}.

\Prop{\label{thm: p1} Proposition \tb{P}} holds for fields $K$ in Table 1 ($h_K=h_F=1$) with $|W_K|=2$ or 4.

\bp
By  Lemmas \ref{thm: quad}, \ref{thm: L2}, \ref{thm: L4} and Remark \ref{thm: sys}, there exists $\wtd{\chi_\mfk{f}}$ satisfying the conditions in Proposition \ref{thm: equiv}. Since $h_K=h_F=1$, $K_\A^\times=K^\times UK^\times_\infty$ and $F^\times_\A=F^\times V F^\times_\infty$. Then $\psi:=\psi_0$ satisfies \tb{a$''$)} through \tb{d$''$)}.
\ep

\Prop{\label{thm: p2} Proposition \tb{P}} holds for fields $K$ in Table 2 ($h_K=h_F=3$) with $|W_K|=2$ or 4.

\bp
By Chebotarev's density theorem, there exists a non-principal prime ideal $\mfk{q}\subset F$ that is inert in $K$, which is equivalent to $\mfk{q}$ being inert in both $H_F$ and $K$. (By \cite{MG}, we can  take $\mfk{q}$ as the prime above (3)). Let $\mfk{Q}\subset K$ be the prime above $\mfk{q}$ and $\lambda\in F$ such that $v_\mfk{q}(\lambda)=1$. Consider $\iota_\mfk{q}\in F_\A^\times$ (resp. $\iota_\mfk{Q}\in K_\A^\times$) which equals to $\lambda$ at the $\mfk{q}$-component (resp. $\mfk{Q}$) and 1 elsewhere. Then since $h_K=h_F=3$, $K_\A^\times =\lan \iota_\mfk{Q}, K^\times U K_\infty^\times\ran $ and $F_\A^\times =\lan \iota_\mfk{q}, F^\times V F^\times_\infty\ran$. Consider $\psi_0$ from Proposition \ref{thm: equiv}. We now extend it to $\psi: K^\times _\A\ra \C^\times $ by setting 
\be
\psi(\iota_\mfk{Q}) := \chi(\iota_\mfk{q})|\iota_\mfk{q}|_{F_\A}^{-1} = \chi(\iota_\mfk{q}) q
\ee
Since $\iota_\mfk{q}^3\in F^\times V F^\times_\infty$ on which \tb{b$''$)} holds, we have $\psi_0(\iota_\mfk{Q}^3)=\psi_0(\iota_\mfk{q}^3)=\chi(\iota_\mfk{q}^3)|\iota_\mfk{q}^3|^{-1}_{F_\A}=\psi(\iota_\mfk{Q})^3$, so this extension is well-defined. By construction, \tb{a$''$)} and \tb{b$''$)} are satisfied. Note that $\psi(\iota_\mfk{Q})\psi(\iota_\mfk{Q})^\rho = q^{2} = |\iota_\mfk{Q}|_{K_\A}^{-1}$. Moreover, we have $(\psi(\iota_\mfk{Q})) = q  =il(f(\iota_\mfk{Q}))$, so \tb{c$''$)} holds. Finally, note that $\mfk{q}^\sigma, \mfk{q}^{\sigma^2}$ are not principal, so $\mfk{q}^\sigma=\mfk{q}^{\delta_1}$ and $\mfk{q}^{\sigma^2}=\mfk{q}^{\delta_2}$ in $C_F$ for some $\delta_1, \delta_2\in \{1, 2\}$. Then as $\mfk{q}\mfk{q}^\sigma\mfk{q}^{\sigma^2}$ is principal, $\delta_1=\delta_2=1$. Since $\iota_\mfk{Q}^\sigma\iota_\mfk{Q}^{-1}\in F^\times V F_\infty^\times$ and $\psi(\iota_\mfk{Q}^\sigma\iota_\mfk{Q}^{-1})=\chi(\iota_\mfk{q}^\sigma\iota_\mfk{q}^{-1})|\iota_\mfk{q}^\sigma\iota_\mfk{q}^{-1}|^{-1}_{F_\A}=1$,  $\psi(\iota_\mfk{Q}^\sigma)=\psi(\iota_\mfk{Q})=\psi(\iota_\mfk{Q})^\sigma$, so \tb{d$''$)} remains valid.
\ep

\Prop{\label{thm: p3} Proposition \tb{P}} holds for fields $K$ in Table 3 ($h_K=4, h_F=1$) with $|W_K|=2$ or 4.

\bp
Let $q\in \Z^+$ be the prime that ramifies in $F_0$. Then by \cite{MG}, we have $q=\mfk{q}_1\mfk{q}_2\mfk{q}_3=\mfk{Q}_1^2\mfk{Q}_2^2\mfk{Q}_3^2$ for primes  $\mfk{q}_j\subset F$ and $\mfk{Q}_j\subset K$. Let $(\eta, \e)=(1+i, 1)$ if $K=\Gamma_1, \Gamma_4$ and $(\sqrt{-q}, 0)$ otherwise. For $1\leq j\leq 3$, note that $K_{\mfk{Q}_j}\cong \Q_q(\eta)$ and $F_{\mfk{q}_j}\cong \Q_q$. Consider $\iota_{\mfk{Q}_1} \in K_\A^\times$ (resp. $\iota_{\mfk{Q}_2}$) that equals to $\eta$ (resp. $\eta^\rho$) at the $\mfk{Q}_1$-component (resp. $\mfk{Q}_2$) and 1 elsewhere. By \cite[Lemma 4]{DC}, $K_\A^\times =\lan \iota_{\mfk{Q}_1}, \iota_{\mfk{Q}_2}, K^\times U K_\infty^\times\ran$ and $\iota_{\mfk{Q}_j}^2\in K^\times UK^\times_\infty$. We extend $\psi_0$ from Proposition \ref{thm: equiv}  to $\psi: K_\A^\times \ra \C^\times $ by setting
\be
\psi(\iota_{\mfk{Q}_1})^\rho=\psi(\iota_{\mfk{Q}_2}):=(-1)^\e\eta
\ee
If $K\neq \Gamma_1, \Gamma_4$, $\psi_0(\iota_{\mfk{Q}_j}^2)=\chi(\iota_{\mfk{Q}_j}^2)|\iota_{\mfk{Q}_j}^2|_{F_\A}^{-1}=\chi_{\mfk{q}_j}(-1)q$, since $N_{K_{\mfk{Q}_j}/F_{\mfk{q}_j}}(\eta)=q$. If $q\neq 2$, since $q\equiv 3$ (mod 4), $\chi_{\mfk{q}_j}(-1)=-1$. This remains true for $q=2$ by the proof of Lemma \ref{thm: L1}. Thus $\psi_0(\iota_{\mfk{Q}_j}^2)=-q=\psi(\iota_{\mfk{Q}_j})^2$. If $K=\Gamma_1, \Gamma_4$, then $\psi_0(\iota_{\mfk{Q}_1}^2)=2\wtd{\chi_{\mfk{f}}}(i, 1, 1)=-2i=\psi(\iota_{\mfk{Q}_1})^2$, as $N_{K_{\mfk{Q}_j}/F_{\mfk{q}_j}}(\eta)=2$ (see Remark \ref{thm: sys}). Similarly, $\psi_0(\iota_{\mfk{Q}_2}^2)=\psi(\iota_{\mfk{Q}_2})^2$. Hence this extension is well-defined. Note that \tb{a$''$)} and \tb{b$''$)} hold as before. Moreover, for $j=1, 2$, $\psi(\iota_{\mfk{Q}_j})\psi(\iota_{\mfk{Q}_j})^\rho=q=|\iota_{\mfk{Q}_j}|_{K_\A}^{-1}$ and $(\psi(\iota_{\mfk{Q}_j}))=(\eta)=\mfk{Q}_1\mfk{Q}_2\mfk{Q}_3=il(f(\iota_{\mfk{Q}_j}))$, so \tb{c$''$)} is satisfied. Finally, $\psi(\iota_{\mfk{Q}_1}^\sigma)=\psi(\iota_{\mfk{Q}_2})=\psi(\iota_{\mfk{Q}_1})^\sigma$. Since $\iota_{\mfk{Q}_1}\iota_{\mfk{Q}_2}\iota_{\mfk{Q}_2}^\sigma\in K^\times U K^\times_\infty$ and in particular $\iota_{\mfk{Q}_1}\iota_{\mfk{Q}_2}\iota_{\mfk{Q}_2}^\sigma\eta^{-1}\in UK^\times_\infty$, 
\be
\psi_0(\iota_{\mfk{Q}_1}\iota_{\mfk{Q}_2}\iota_{\mfk{Q}_2}^\sigma)=\wtd{\chi_{\mfk{f}}}(1, \frac{\eta^\rho}{\eta}, 1) (f(\eta^{-1})_{\infty 1})^{-1}=\wtd{\chi_{\mfk{f}}}(1, \frac{\eta^\rho}{\eta}, 1)\eta\eta^{\sigma^4}\eta^{\sigma^5}=\wtd{\chi_{\mfk{f}}}(1, \frac{\eta^\rho}{\eta}, 1)\eta^2\eta^\rho
\ee
If $K\neq \Gamma_1, \Gamma_4$, then $\eta^\rho\eta^{-1} = -1$, so $\psi(\iota^\sigma_{\mfk{Q}_2})=-\eta=\psi(\iota_{\mfk{Q}_2})^\sigma$. If $K=\Gamma_1, \Gamma_4$, $\eta^\rho\eta^{-1}=-i$. As $\wtd{\chi_{\mfk{f}}}(1, -i, 1)=i$, we have $\psi(\iota_{\mfk{Q}_2}^\sigma)=i\eta=(-\eta)^\sigma=\psi(\iota_{\mfk{Q}_2})^\sigma$. Hence \tb{d$''$)} is valid.
\ep

\Prop{\label{thm: p4} Proposition \tb{P}} holds for $K=\Omega_1$ in Table 4 ($h_K=12, h_F=3$).
\bp
Suppose $2=\mfk{q}_1\mfk{q}_2\mfk{q}_3=\mfk{Q}_1^2\mfk{Q}_2^2\mfk{Q}_3^2$ for primes $\mfk{q}_j\subset F$, $\mfk{Q}_j\subset K$. Let $\mfk{r}\subset F$ be a non-principal prime  that is inert in $K$ and $\mfk{R}\subset K$ the prime above it. Let $\eta=1+i$ and $\lambda\in F$ such that $v_\mfk{r}(\lambda)=1$. Consider $\iota_{\mfk{Q}_1}\in K^\times_\A$ (resp. $\iota_{\mfk{Q}_2}, \iota_\mfk{R}$) that equals to $\eta$ (resp. $\eta^\rho, \lambda$) at the $\mfk{Q}_1$-component (resp. $\mfk{Q}_2, \mfk{R}$) and 1 elsewhere. Note that $K^\times_\A=\lan \iota_{\mfk{Q}_1}, \iota_{\mfk{Q}_2}, \iota_{\mfk{R}}, K^\times U K^\times_\infty\ran$. We extend $\psi_0$ from Proposition \ref{thm: equiv} by setting $\psi(\iota_{\mfk{Q}_1})^\rho=\psi(\iota_{\mfk{Q}_2}):=-\eta$ and $\psi(\iota_\mfk{R}):=\chi(\iota_\mfk{r})|\iota_\mfk{r}|_{F_\A}^{-1}$, where $\iota_\mfk{r}$ is $\iota_{\mfk{R}}$ considered as an idele in $F_\A^\times$. Then by the arguments in Propositions \ref{thm: p2} and  \ref{thm: p3}, we have $\psi$ is well-defined and satisfies \tb{a$''$)} through \tb{d$''$)}.
\ep

\Prop{\label{thm: cond} If $|W_K|\neq 6$, the $\Q$-models of $(A, \theta)$ have conductors $\cond(A_{/\Q})$ as listed in Tables \ref{tab: T1} through \ref{tab: T4}. Moreover, the conductors of $\psi$ are denoted $\cond(\psi)$.

\bp
Let $\cond(A_{/K})$ be the conductor of $A$ over $K$. By \cite[Theorem 12]{ST} and \cite[Section 3, Corollary (b)]{JM}, we have
\begin{align*}
\cond(A_{/K}) &=\cond(\psi)^6\\
N_{K/\Q}(\cond(A_{/K})) d_K^6 &=\cond(A_{/\Q})^6
\end{align*}
Then we can compute $\cond(A_{/\Q})$ from $d_K$ and $\cond(\psi)$, the latter being immediate from our construction. In the tables, $\mfk{P}, \mfk{Q}$ are primes in $K$ above the first and second primes in the factorization of $|d_K|$.
\ep

\vspace{-\len}
Now we show that the desired Hecke characters don't exists for fields with $|W_K|=6$.

\Rem{\label{thm: wk6} For $K=A_4, A_7, A_{15}, A_{17}$, the prime ideal (3) is inert in $F$. For $B_1, B_2, B_5$, it ramifies, while for $\Gamma_5, \Gamma_6$, it splits completely. Let $\mfk{p}\subset F$ be a prime above $(3)$. Then $\mfk{p}\mcl{O}_K=\mfk{P}^2$ for some prime $\mfk{P}\subset K$, and $\mfk{f}=\mfk{p}$ or $\mfk{p}_1\mfk{p}_2\mfk{p}_3$. Let $q\ (\neq p)$ be the rational prime dividing the conductor of $F$. Then $q\mcl{O}_F=\mfk{q}^3$ for a prime $\mfk{q}\subset F$ and $\mfk{q}$ splits in $K$. (These follow from direct computation with Magma \cite{MG}).}

\Lem{\label{thm: kum} For $K=A_4, A_7, A_{15}, A_{17}, B_1, B_2, B_5$ and $\mfk{P}$ as in Remark \ref{thm: wk6}, if $\wtd{\chi_\mfk{P}}: \mcl{O}_{K_\mfk{P}}^\times \ra W_K$ is a homomorphism with $\wtd{\chi_\mfk{P}}(x^\sigma)=\wtd{\chi_\mfk{P}}(x)^\sigma$, then $\wtd{\chi_\mfk{P}}(\omega)=1$.}

\bp
Note that $K_\mfk{P}$ is a cubic extension of $\Q_3(\omega)$. By Kummer theory, there exists a size three subgroup $\lan \alpha \ran$ of $\Q_3(\omega)^\times/\Q_3(\omega)^{\times 3}$ generated by $\alpha$ such that $K_\mfk{P}=\Q_3(\omega, \sqrt[3]{\alpha})$. Since $K_\mfk{P}$ is cyclic over $\Q_3$,  $\lan \alpha\ran$ is stable under the action of $\gal(\Q_3(\omega)/\Q_3)=\lan \rho\ran$, and $\rho$ acts nontrivially on $\alpha$, so $\alpha^\rho=\alpha^2$ in $\Q_3(\omega)^\times/\Q_3(\omega)^{\times 3}$. This implies $\alpha\in \mcl{O}_{\Q_3(\omega)}^\times$, so $\sqrt[3]{\alpha}\in \mcl{O}_{K_\mfk{P}}^\times$. Since $\wtd{\chi_\mfk{P}}(\omega\sqrt[3]{\alpha})=\wtd{\chi_\mfk{P}}(\sqrt[3]{\alpha}^{\sigma^2})=\wtd{\chi_\mfk{P}}(\sqrt[3]{\alpha})^{\sigma^2}=\wtd{\chi_\mfk{P}}(\sqrt[3]{\alpha})$, $\wtd{\chi_\mfk{P}}(\omega)=1$
\ep

\Prop{\label{thm: no} Proposition \tb{P} holds for fields $K$ in Tables 1 through 4 with $|W_K|=6$.}

\bp
Suppose there is a Hecke character $\psi$ of $K^\times_\A$ satisfying \tb{a$''$)} through \tb{d$''$)}. Let $\mfk{p}\subset F$ be a prime  and $\mfk{P}\subset K$ (resp. $p\in \Z^+$) the prime above (resp. below) it. Let $\psi_\mfk{P}$ be the character on $K_\mfk{P}^\times$ induced by $\psi$. Now we show that $\prod_u\psi_u(\omega)=1$, where $u$ runs through all the finite places of $K$. There are three cases:

\tb{i)} $\mfk{p}\mcl{O}_K=\mfk{P}_1\mfk{P_2}$. Note that $K_{\mfk{P}_1}=F_\mfk{p}=K_{\mfk{P}_2}$. Since $\mfk{p}\nmid \mfk{f}$,  we have $\psi_{\mfk{P}_1}(\omega)\psi_{\mfk{P}_2}(\omega)=\chi_\mfk{p}(\omega)=1$ by \tb{b$''$)}.

\tb{ii)} $\mfk{p}\mcl{O}_K=\mfk{P}$. If $p\mcl{O}_K=\mfk{P}_1\mfk{P}_2\mfk{P}_3$, suppose $\sigma^2\mfk{P}_1=\sigma\mfk{P}_2=\mfk{P}_3$ and $\psi_{\mfk{P}_1}(\omega^2)=\lambda\in\{ 1, \omega, \omega^2\}$. Then by \tb{d$''$)}, $\psi_{\mfk{P}_2}(\omega)=\psi_{\mfk{P}_1}(\omega^2)^\sigma=\lambda^\rho$ and $\psi_{\mfk{P}_3}(\omega^2)=\lambda^{\sigma^2}=\lambda$, so $\prod_{j=1}^3\psi_{\mfk{P}_j}(\omega)=\lambda^4\lambda^\rho=1$. Otherwise, since $\mfk{q}$ splits in $K$ by Remark \ref{thm: wk6}, we have $p\mcl{O}_K=\mfk{P}$. Note that the residual field of $\mcl{O}_{K_\mfk{P}}$ has size $p^6$. Since $9|p^6-1$, by Hensel's Lemma, there exists $\zeta_9\in \mcl{O}_{K_\mfk{P}}^\times$, so $\psi_\mfk{P}(\omega)= 1$.

\tb{iii)} $\mfk{p}\mcl{O}_K=\mfk{P}^2$. Then $p=3$. If $K\neq\Gamma_5, \Gamma_6$, $\forall \lambda\in \mcl{O}_{K_\mfk{P}}^\times$, we have $\psi_\mfk{P}(\lambda)\psi_\mfk{P}(\lambda)^\rho=1$ by \tb{c$''$)}, so $\psi_\mfk{P}(\lambda)\in W_K$.
Then $\psi_\mfk{P}(\omega)=1$ by Lemma \ref{thm: kum}. If $K=\Gamma_5, \Gamma_6$, suppose $3\mcl{O}_K=\mfk{P}_1^2\mfk{P}_2^2\mfk{P}_3^2$. Then the argument in the first half of \tb{ii)} shows that $\prod_{j=1}^3\psi_{\mfk{P}_j}(\omega)=1$.

Thus $\prod_u\psi_u(\omega)= 1$. Let $\iota\in K_\infty^\times$ be the idele equal to the $\infty$-component of $\omega \in K^\times_\A$. Then for $\omega\in K^\times$, $\psi(\omega)=\psi(\iota)\prod_u\psi_u(\omega)=(f(\iota)_{\infty 1})^{-1}=(\omega\omega^{\sigma^4}\omega^{\sigma^5})^{-1}=\omega^2\neq 1$. Contradiction! Hence no such $\psi$ exists.
\ep

\Prop{\label{thm: last} For fields $K$ in Tables 1 through 4 with $|W_K|=6$, $(A, \theta)$ has a $F$-model.}
\bp
Take $D=D'=F=k_0$ in Theorem \ref{thm: S1}. It suffices to construct a Hecke character $\psi: K^\times_\A\ra \C^\times$ that satisfies \tb{a$''$)}, \tb{b$''$)}, \tb{c$''$)} and Galois equivariance by $\rho$. The corresponding local character $\wtd{\chi_\mfk{f}}$ can be constructed explicitly in a fashion similar to the proof of Lemmas \ref{thm: L1} and \ref{thm: L3}.
\ep

\section{Future Directions}
An immediate direction for future research is to verify if the analogue to $(\ast)$ from the abstract holds for four dimensional CM abelian varieties with rational fields of moduli. This will probably involve an even longer list of degree eight CM fields than Tables \ref{tab: T1} through \ref{tab: T4}. As we can see from the proof above, a large $|W_K|$ is the major obstruction to the existence of our desired Hecke charcters, which is made precise by an obstruction class in $H^2$ in \cite[Proposition 1]{HY}. Thus, it is unlikely that these abelian varieties all have $\Q$-models.

Alternatively, we can fix an integer $n\geq 4$ and ask if there is a finite list of degree $2n$ CM fields from which $n$-dimensional CM abelian varieties with fields of moduli $=\Q$ could arise. For $n=2, 3$, the approaches by Murabayashi \cite{M1} and Chun \cite{DC} start with the identification of the relative class numbers of $K/F$. Hence an affirmative answer may hinge on the classification of extensions of a fixed relative class number as in \cite{PK}.

Finally, recall that Jacobi sum characters (discovered by Weil and formalized by Anderson) are  a class of Hecke charcters systematically constructed for all abelian fields, with known connection to Fermat hypersurfaces.
Since our $K$ is abelian, it is of interest to figure out which  Hecke characters constructed for Proposition \tb{P} arise from Jacobi sums, which could improve our understanding of the corresponding abelian varieties.

\section*{Appendix}

\begin{small}
\vspace{-12pt}

\begin{table}[H]
\caption{Degree six CM fields $K$ with $h_K=1$ and $h_{F}=1$}
\label{tab: T1}
\begin{center}
\begin{tabular}{cccccccccc}
\hline
$K$ & $f_K$ & $f_{F}$ & $f_{F_0}$ &  $P_F$ & $|d_K|$ & $\cond(\psi)$  & $\cond(A_{/\Q})$ \\
\hline
$\Q(\alpha_1, \sqrt{-7})$ & 7 & 7 & 7 &  $x^3 - 21x + 7$ & $7^5$ & $\mfk{P}^2$ &  $7^7$\\
$\Q(\alpha_2, \sqrt{-3})$ & 9 & 9 & 3 &  $x^3-3x-1$ & $3^9$ & $\mfk{P}^4$ & $3^{13}$\\
$\Q(\alpha_3, \sqrt{-19})$ & 19 & 19 & 19 &  $x^3 - 57x + 133$ & $19^5$ & $\mfk{P}$ & $19^6$\\
$\Q(\alpha_4, \sqrt{-3})$ & 21 & 7 & 3 &  $x^3 - 21x + 7$ & $3^3\cdot 7^4$ & &\\
$\Q(\alpha_5, \sqrt{-1})$ & 28 & 7 & 4 &  $x^3 -21x+7$ & $2^6\cdot 7^4$ & $\mfk{P}^3$ & $2^{15}\cdot 7^4$\\
$\Q(\alpha_6, \sqrt{-1})$ & 36 & 9 & 4 &  $x^3-3x-1$ & $2^6\cdot 3^8 $ & $\mfk{P}^3$ & $2^{15}\cdot 3^8$\\
$\Q(\alpha_7, \sqrt{-3})$ & 39 & 13 & 3 &  $x^3-39x -65$ & $3^3\cdot 13^4$ & &\\
$\Q(\alpha_8, \sqrt{-43})$ & 43 & 43 & 43 &  $x^3-129x-344$ & $43^5$ & $\mfk{P}$ & $43^6$\\
$\Q(\alpha_9, \sqrt{-2})$ & 56 & 7 & 8 &  $x^3-21x+7$ & $2^9\cdot 7^4$ & $\mfk{P}^5$ & $2^{24}\cdot 7^4$\\
$\Q(\alpha_{10}, \sqrt{-7})$ & 63 & 9 & 7 &  $x^3-3x-1$ & $3^8\cdot 7^3$ & $\mfk{Q}$ & $3^8\cdot 7^6$\\
$\Q(\alpha_{11}, \sqrt{-67})$ & 67 & 67 & 67 &  $x^3-201x-335$ & $67^5$ & $\mfk{P}$ & $67^6$\\
$\Q(\alpha_{12}, \sqrt{-1})$ & 76 & 19 & 4 &  $x^3-57x+133$ & $2^6\cdot 19^4$ & $\mfk{P}^3$ & $2^{15}\cdot 19^4$\\
$\Q(\alpha_{13}, \sqrt{-11})$ & 77 & 7 & 11 &  $x^3-21x+7$ & $7^4\cdot 11^3$ & $\mfk{Q}$ & $7^4\cdot 11^6$\\
$\Q(\alpha_{14}, \sqrt{-7})$ & 91 & 13 & 7 &  $x^3-39x-65$ & $7^3\cdot 13^4$ & $\mfk{P}$ & $7^6\cdot 13^4$\\
$\Q(\alpha_{15}, \sqrt{-3})$ & 93 & 31 & 3 &  $x^3-93x+124$ & $3^3\cdot 31^4$ & &\\
$\Q(\alpha_{16}, \sqrt{-2})$ & 104 & 13 & 8 &  $x^3-39x-65$ & $2^9\cdot 13^4$ & $\mfk{P}^5$ & $2^{24}\cdot 13^4$\\
$\Q(\alpha_{17}, \sqrt{-3})$ & 129 & 43 & 3 &  $x^3-129x-344$ & $3^3\cdot 43^4$ & &\\
\hline
\end{tabular}
\end{center}
\end{table}

\vspace{-12pt}

\begin{table}[H]
\caption{Degree six CM fields $K$ with $h_K=3$ and $h_{F}=3$}
\label{tab: T2}
\begin{center}
\begin{tabular}{ccccccccc}
\hline
$K$ & $f_K$ & $f_{F}$ & $f_{F_0}$ &  $P_F$ & $|d_K|$ &$\cond(\psi)$ &  $\cond(A_{/\Q})$\\
\hline
$\Q(\beta_1, \sqrt{-3})$ & 63 & 63 & 3 &  $x^3 - 21x -28$ & $3^9\cdot 7^4$\\
$\Q(\beta_2, \sqrt{-3})$ & 63 & 63 & 3 &  $x^3 -21x+35$ & $3^9\cdot 7^4$\\
$\Q(\beta_3, \sqrt{-7})$ & 63 & 63 & 7 &  $x^3 -21x +35$ & $3^8\cdot 7^5$ & $\mfk{Q}$ & $3^8\cdot 7^6$\\
$\Q(\beta_4, \sqrt{-7})$ & 91 & 91 & 7 &  $x^3 - x^2 -30x - 27$ & $7^5\cdot 13^4$ & $\mfk{P}$ & $7^6\cdot 13^4$\\
\hline
\end{tabular}
\end{center}
\end{table}

\begin{table}[H]
\begin{center}
\begin{tabular}{cccccccccc}
\hline
$K$ & $f_K$ & $f_{F}$ & $f_{F_0}$ &  $P_F$ & $|d_K|$ & $\cond(\psi)$ & $\cond(A_{/\Q})$\\
\hline
$\Q(\beta_5, \sqrt{-3})$ & 117 & 117 & 3 &  $x^3 -39x+26$ & $3^9\cdot 13^4$\\
$\Q(\beta_6, \sqrt{-7})$ & 133 & 133 & 7 &  $x^3-x^2-44x-69$ & $7^5\cdot 19^4$ & $\mfk{P}$ & $7^6\cdot 19^4$\\
$\Q(\beta_7, \sqrt{-19})$ & 171 & 171 & 19 &  $x^3-57x-152$ & $3^8\cdot 19^5$ & $\mfk{Q}$ & $3^8\cdot 19^6$\\
$\Q(\beta_8, \sqrt{-7})$ & 217 & 217 & 7 &  $x^3-x^2-72x+225$ & $7^5\cdot 31^4$ & $\mfk{P}$ & $7^6\cdot 31^4$\\
$\Q(\beta_9, \sqrt{-19})$ & 247 & 247 & 19 &  $x^3-x^2-82x+64$ & $13^4\cdot 19^5$ & $\mfk{Q}$ & $13^4\cdot 19^6$\\
\hline
\end{tabular}
\end{center}
\end{table}

\vspace{-12pt}

\begin{table}[H]
\caption{Degree six CM fields $K$ with $h_K=4$ and $h_{F}=1$}
\label{tab: T3}
\begin{center}
\begin{tabular}{cccccccccc}
\hline
$K$ & $f_K$ & $f_{F}$ & $f_{F_0}$ & $P_F$ & $|d_K|$ &$\cond(\psi)$ &  $\cond(A_{/\Q})$ \\
\hline
$\Q(\gamma_1, \sqrt{-1})$ & 124 & 31 & 4 &  $x^3 - 93x + 124$ & $2^6\cdot 31^4$ & $\mfk{P}_1^3\mfk{P}_2^3\mfk{P}_3^3$ & $2^{15}\cdot 31^4$\\
$\Q(\gamma_2, \sqrt{-7})$ & 133 & 19 & 7 &  $x^3 - 57x + 133$ & $7^3\cdot 19^4$ & $\mfk{P}_1\mfk{P}_2\mfk{P}_3$ & $7^6\cdot 19^4$ \\
$\Q(\gamma_3, \sqrt{-19})$ & 171 & 9 & 19 &  $x^3 - 3x-1$ & $3^8\cdot 19^3$ & $\mfk{Q}_1\mfk{Q}_2\mfk{Q}_3$ & $3^8\cdot 19^6$\\
$\Q(\gamma_4, \sqrt{-1})$ & 172 & 43 & 4 &  $x^3 - 129x - 344$ & $2^6\cdot 43^4$ & $\mfk{P}_1^3\mfk{P}_2^3\mfk{P}_3^3$ & $2^{15}\cdot 43^4$\\ 
$\Q(\gamma_5, \sqrt{-3})$ & 183 & 61 & 3 &  $x^3 - 183x + 61$ & $3^3\cdot 61^4$\\
$\Q(\gamma_6, \sqrt{-3})$ & 201 & 67 & 3 &  $x^3 - 201x - 335$ & $3^3\cdot 67^4$\\
$\Q(\gamma_7, \sqrt{-11})$ & 209 & 19 & 11 &  $x^3 -57x + 133$ & $11^3\cdot 19^4$ & $\mfk{P}_1\mfk{P}_2\mfk{P}_3$ & $11^6\cdot 19^4$\\
$\Q(\gamma_8, \sqrt{-2})$ & 248 & 31 & 8 &  $x^3 -93x + 124$ & $2^9\cdot 31^4$ & $\mfk{P}_1^5\mfk{P}_2^5\mfk{P}_3^5$ & $2^{24}\cdot 31^4$\\
$\Q(\gamma_9, \sqrt{-11})$ & 473 & 43 & 11 &  $x^3 - 129x - 344$ & $11^3\cdot 43^4$ &$\mfk{P}_1\mfk{P}_2\mfk{P}_3$ & $11^6\cdot 43^4$\\
$\Q(\gamma_{10}, \sqrt{-7})$ & 511 & 73 & 7 &  $x^3 -219x + 511$ & $7^3\cdot 73^4$  &$\mfk{P}_1\mfk{P}_2\mfk{P}_3$ & $7^6\cdot 73^4$\\
\hline
\end{tabular}
\end{center}
\end{table}

\vspace{-12pt}

\begin{table}[H]
\caption{Degree six CM fields $K$ with $h_K=12$ and $h_{F}=3$}
\label{tab: T4}

\begin{center}
\begin{tabular}{ccccccccc}
\hline
$K$ & $f_K$ & $f_{F}$ & $f_{F_0}$ & $P_F$ & $|d_K|$ &$\cond(\psi)$ &  $\cond(A_{/\Q})$ \\
\hline
$\Q(\omega_1, \sqrt{-1})$ & 252 & 63 & 4 &  $x^3 - 21x -28$ & $2^6\cdot 3^8 \cdot 7^4$& $\mfk{P}_1^3\mfk{P}_2^3\mfk{P}_3^3$ & $2^{15}\cdot 3^8\cdot 7^4$\\
\hline
\end{tabular}
\end{center}
\end{table}

\end{small}

\noindent
Department of Mathematics, California Institute of Technology, Pasadena, CA 91125\\
Email address: zshang@caltech.edu

\end{document}